\documentclass{amsart}
\usepackage{amsfonts}

\usepackage{amscd}

\newtheorem{theorem}{Theorem}[section]
\newtheorem{corollary}[theorem]{Corollary}

\newtheorem{lemma}[theorem]{Lemma}

\theoremstyle{definition}

\theoremstyle{remark}

\begin{document}
\title{Cohomology and Obstructions III: A variational form of the generalized Hodge
conjecture on $K$-trivial threefolds}
\author{Herb Clemens}
\address{Mathematics Department, University of Utah}
\email{clemens@math.utah.edu}
\date{February, 2002}
\maketitle

\begin{abstract}
This paper studies the Hilbert scheme of a curve on a complete-intersection $%
K$-trivial threefold, in the case in which the curve is unobstructed in the
ambient variety in which the threefold lives. The basic result is that the
obstruction theory of the curve is completely determined by the
scheme-theoretic Abel-Jacobi mapping. Several applications of this fact are
given.
\end{abstract}

\section{Introduction}

This paper gives a proof of a variational version of the generalized Hodge
conjecture in a very special case, namely for certain families of curves on
certain complex projective threefolds $X_{0}$ with trivial canonical bundle.
The (smooth) threefold $X_{0}$ must be the zero-scheme of a regular section
of a vector bundle on a complex manifold $P$ and a maximal abelian
subvariety of the intermediate Jacobian $J\left( X_{0}\right) $ must deform
over some deformation $X/X^{\prime }$ of $X_{0}$ as a sub-manifold of $P$.
Finally, over a polydisk $Y_{0}^{\prime }$, we must be given a family of
(smooth) curves
\begin{equation*}
p:Y_{0}\rightarrow Y_{0}^{\prime }
\end{equation*}
in $X_{0}$ whose fibers are (strongly) unobstructed in $P$, that is,
\begin{equation*}
R^{1}p_{*}N_{Y_{0}\backslash Y_{0}^{\prime }\times P}.
\end{equation*}
The conclusion is then that $Y_{0}/Y_{0}^{\prime }$ is the specialization of
some family of curves on the generic fiber of $X/X^{\prime }$. This result
should be inerpreted as generalizing the fact that a rigid curve on $X_{0}$
is always the specialization of a curve on the general fiber of $X/X^{\prime
}$ (because the relative dimension of the Hilbert scheme is bounded below by
the Euler characteristic of the normal bundle which, for $K$-trivial
threefolds, is alway zero).

A previous version of this paper appeared with the title ``Cohomology and
Obstructions II: Curves on Calabi-Yau threefolds.'' C. Voisin subsequently
pointed out that a key assertion in the proof of the main theorem, namely
that an Abel-Jacobi map over a scheme takes values in an abelian variety
defined over that scheme, was unsubstantiated. In developing a proof that
claim in the case of intermediate Jacobians of $K$-trivial threefolds, the
author was led to develop an algebro-geometric version of the physicists'
understanding of Hilbert schemes of Calabi-Yau threefolds as gradient
schemes. That study developed a life of its own, resulting in a separate
paper \cite{C2}, some of the results of which are necessary to correct the
above-mentioned gap. To keep things in logical order, this manuscript has
been renamed, and carries the label ``Cohomology and Obstructions III.''

\subsection{The tool: Gauss-Manin connection as obstruction}

This paper represents a concrete application of Theorem 13.1 of \cite{C1}
and of the realization of the relative Hilbert scheme as a relative gradient
varient, which is the main result of \cite{C2}. We first recall Theorem 13.1
of \cite{C1}. In what follows, 
\begin{equation*}
X/X^{\prime }
\end{equation*}
will always be a deformation of a complex projective manifold $X_{0}$ where $%
X^{\prime }$ is a complex polydisk with distinguished basepoint $0$. (Later
in this paper, $X_{0}$ will be a $K$-trivial threefold, like a Calabi-Yau
threefold or an abelian variety.) For any base extension 
\begin{equation*}
A\rightarrow X^{\prime }
\end{equation*}
we let 
\begin{equation*}
X_{A}:=X\times _{X^{\prime }}A.
\end{equation*}

\begin{theorem}
\label{k}Let $\tilde{\Delta }$ denote a small disk around $0$ in $\Bbb{C}$
with parameter $\tilde{t}$. Let $X_{\tilde{\Delta }}/\tilde{\Delta }$ be a
deformation of the complex projective manifold $X_{0}$. Suppose 
\begin{equation*}
p^{\prime }:W_{0}\rightarrow W_{0}^{\prime }
\end{equation*}
is a proper family of submanifolds of $X_{0}$ of fiber dimension $q$ over a
smooth (not necessarily compact) base $W_{0}^{\prime }$. Suppose further
that the family $W_{0}/$ $W_{0}^{\prime }$ deforms with $X_{0}$ to a family $%
W_{n}/$ $W_{n}^{\prime }$ over the Artinian scheme $\tilde{\Delta }%
_{n}\subseteq \tilde{\Delta }$ associated to the ideal $\left\{ \tilde{t}%
^{n+1}\right\} $ and that 
\begin{equation*}
\omega \in H^{p+q+1,q-1}\left( X_{\tilde{\Delta }}/\tilde{\Delta }\right) 
\end{equation*}
lies in the kernel of the composition 
\begin{equation*}
H^{p+q+1,q-1}\left( X_{\tilde{\Delta }}/\tilde{\Delta }\right) \overset{%
\nabla }{\longrightarrow }H^{p+q,q}\left( X_{\tilde{\Delta }}/\tilde{\Delta }%
\right) \overset{\left( pull-back\right) }{\longrightarrow }%
R^{q}p_{*}^{\prime }\left( \Omega _{W_{n}/\tilde{\Delta }_{n}}^{p+q}\right)
\rightarrow \Omega _{W_{n}^{\prime }}^{p}
\end{equation*}
induced by the Gauss-Manin connection $\nabla $ and integration over the
fiber. Let 
\begin{equation}
\varsigma \in R^{1}p_{*}^{\prime }\left( N_{W_{0}\backslash W_{0}^{\prime }%
\times X_{0}}\right)   \label{genobst}
\end{equation}
be the obstruction class measuring extendability of $W_{n}/$ $W_{n}^{\prime }
$ to a family $W_{n+1}/$ $W_{n+1}^{\prime }.$ Then, for $\omega _{0}=\left.
\omega \right| _{X_{0}},$%
\begin{equation}
\left. \left\langle \left. \varsigma \right| \omega _{0}\right\rangle
\right| _{W_{0}}\in R^{q}p_{*}^{\prime }\left( \Omega _{W_{0}}^{p+q}\right)
\otimes \frac{_{\left\{ \tilde{t}^{n+1}\right\} }}{_{\left\{ \tilde{t}%
^{n+2}\right\} }}  \label{BF-1}
\end{equation}
goes to zero in 
\begin{equation*}
\Omega _{W_{0}^{\prime }}^{p}\otimes \frac{_{\left\{ \tilde{t}^{n+1}\right\}
}}{_{\left\{ \tilde{t}^{n+2}\right\} }}
\end{equation*}
under the map 
\begin{equation}
\int\nolimits_{W_{0}/W_{0}^{\prime }}:R^{q}p_{*}^{\prime }\left( \Omega %
_{W_{0}}^{p+q}\right) \rightarrow \Omega _{W_{0}^{\prime }}^{p}  \label{BF-2}
\end{equation}
induced by the Leray spectral sequence. (This last map is commonly called 
\textit{integration over the fiber}.)
\end{theorem}

\subsection{Abel-Jacobi mapping\label{AJ}}

Suppose again that 
\begin{equation*}
X_{\tilde{\Delta }}/\tilde{\Delta }
\end{equation*}
is a deformation of $X_{0}$ over the $\tilde{t}$-disk. Let $W_{0}^{\prime }$
be a complex polydisk. Suppose we have a diagram 
\begin{equation*}
\begin{array}{ccc}
W_{0} & \hookrightarrow & X_{0} \\ 
\downarrow p_{0} &  &  \\ 
W_{0}^{\prime } &  & 
\end{array}
\end{equation*}
where $p_{0}$ is proper with fibers 
\begin{equation*}
\left\{ Z_{w}\right\} _{w\in W_{0}^{\prime }}
\end{equation*}
which are smooth and irreducible of dimension $q$ and suppose, for
simplicity, that the total space $W_{0}$ of the family is embedded in $X_{0} 
$. Let 
\begin{equation*}
\tilde{\Delta }_{n}:=Spec\frac{\Bbb{C}\left[ \left[ \tilde{t}\right] \right] 
}{\left\{ \tilde{t}^{n+1}\right\} }.
\end{equation*}
Suppose $p_{0}$ extends to 
\begin{equation}
p_{n}:W_{n}\rightarrow W_{n}^{\prime }=W_{0}^{\prime }\times \tilde{\Delta }%
_{n}  \label{nf1}
\end{equation}
a proper and smooth family with fiber dimension $q$ embedded in $X_{n}$, the
part of $X_{\tilde{\Delta }}$ over $\tilde{\Delta }_{n}$ via the commutative
diagram 
\begin{equation*}
\begin{array}{ccc}
W_{n} & \hookrightarrow & X_{n} \\ 
\downarrow p_{n} &  & \downarrow \\ 
W_{n}^{\prime } & \rightarrow & \tilde{\Delta }_{n}
\end{array}
.
\end{equation*}

Let 
\begin{equation*}
\mathcal{J}_{q}\left( X_{\tilde{\Delta }}/\tilde{\Delta }\right) =\frac{%
\left( H^{0}\left( \Omega _{X_{\tilde{\Delta }}/\tilde{\Delta }%
}^{2q+1}\right) +\ldots +H^{q+1}\left( \Omega _{X_{\tilde{\Delta }}/\tilde{%
\Delta }}^{q}\right) \right) ^{\vee }}{H_{2q+1}\left( X_{\tilde{\Delta }}/%
\tilde{\Delta };\Bbb{Z}\right) }
\end{equation*}
be the relative $q$-th intermediate Jacobian of $X/\Delta $. Then we have an
Abel-Jacobi mapping 
\begin{equation*}
\begin{array}{c}
\varphi :W_{0}^{\prime }\times \tilde{\Delta }_{n}\rightarrow \mathcal{J}%
_{q}\left( X_{\tilde{\Delta }}/\tilde{\Delta }\right) \\ 
\left( w,\tilde{t}\right) \mapsto \int\nolimits_{Z_{0}}^{Z_{w}}
\end{array}
\end{equation*}
defined as a morphism of analytic schemes. Further suppose 
\begin{equation}
A/\tilde{\Delta }\subseteq \mathcal{J}_{q}\left( X_{\tilde{\Delta }}/\tilde{%
\Delta }\right)  \label{split}
\end{equation}
is an abelian subvariety orthogonal to $H^{0}\left( \Omega _{X_{\tilde{%
\Delta }}/\tilde{\Delta }}^{2q+1}\right) +\ldots +H^{q-1}\left( \Omega _{X_{%
\tilde{\Delta }}/\tilde{\Delta }}^{q+2}\right) $ such that 
\begin{equation*}
\varphi \left( W_{n}^{\prime }\right) \subseteq A_{n}
\end{equation*}
the abelian fiber over $\tilde{t}=0$. Then $\left( \ref{split}\right) $
induces a splitting of 
\begin{equation*}
\mathcal{J}_{q}\left( X_{\tilde{\Delta }}/\tilde{\Delta }\right) =\frac{%
\left( H^{0}\left( \Omega _{X_{\tilde{\Delta }}/\tilde{\Delta }%
}^{2q+1}\right) +\ldots +H^{q+1}\left( \Omega _{X_{\tilde{\Delta }}/\tilde{%
\Delta }}^{q}\right) \right) ^{\vee }}{H_{2q+1}\left( X_{\tilde{\Delta }}/%
\tilde{\Delta };\Bbb{Z}\right) }
\end{equation*}
as follows. By Hodge theory there exists a lattice 
\begin{equation*}
L=L_{1}\oplus L_{2}
\end{equation*}
in $H_{2q+1}\left( X_{\tilde{\Delta }}/\Delta ;\Bbb{Q}\right) $ such that $%
H_{2q+1}\left( X_{\tilde{\Delta }}/\Delta ;\Bbb{Z}\right) $ is of finite
index and a direct sum decomposition 
\begin{equation*}
H^{0}\left( \Omega _{X_{\tilde{\Delta }}/\tilde{\Delta }}^{2q+1}\right)
+\ldots +H^{q+1}\left( \Omega _{X_{\tilde{\Delta }}/\tilde{\Delta }%
}^{q}\right) =\left( V_{1}/\tilde{\Delta }\right) \oplus \left( V_{2}/\tilde{%
\Delta }\right)
\end{equation*}
such that 
\begin{equation*}
A/\tilde{\Delta }=\frac{V_{1}^{\vee }}{L_{1}}
\end{equation*}
and the above decomposition induces an isogeny 
\begin{equation*}
\mathcal{J}_{q}\left( X_{\tilde{\Delta }}/\tilde{\Delta }\right) \rightarrow
\left( A\oplus A^{\bot }\right) /\tilde{\Delta },
\end{equation*}
where 
\begin{equation*}
A^{\bot }=\frac{V_{2}^{\vee }}{L_{2}}.
\end{equation*}
Furthermore the action of the Gauss-Manin connection respects the direct sum
decompostion. Thus, since 
\begin{equation*}
H^{q+2,q-1}\left( X_{\tilde{\Delta }}/\tilde{\Delta }\right) \subseteq V_{2}/%
\tilde{\Delta },
\end{equation*}
differentiating via the Gauss-Manin connection, we have 
\begin{equation*}
\frac{\nabla H^{q+2,q-1}\left( X_{\tilde{\Delta }}/\tilde{\Delta }\right) }{%
\partial \tilde{t}}\subseteq V_{2}/\tilde{\Delta }.
\end{equation*}
So from Theorem 13.1 of \cite{C1} we have:

\begin{corollary}
\label{goodcor}If, under the Abel-Jacobi mapping 
\begin{equation*}
\varphi \left( W_{n}^{\prime }\right) \subseteq A_{n},
\end{equation*}
the mapping 
\begin{equation*}
R^{1}p_{*}^{\prime }\left( \Omega _{W_{0}}^{3}\right) \rightarrow \Omega %
_{W_{0}^{\prime }}^{1}
\end{equation*}
given by $\left( \ref{BF-2}\right) $ is zero.
\end{corollary}

\subsection{\label{set}The setting: Curves on $K$-trivial threefolds}

Let $\Delta $ denote the $t$-disk and $\Delta _{n}\subseteq \Delta $ the
subscheme associated to the ideal $\left\{ t^{n+1}\right\} .$ Our special
situation is one in which there is a fixed projective manifold $P$, a vector
bundle 
\begin{equation*}
\frak{E}\rightarrow P,
\end{equation*}
and 
\begin{equation*}
X_{\Delta }\subseteq P\times \Delta
\end{equation*}
is a family of smooth threefolds $\left\{ X_{t}\right\} _{t\in \Delta }$
given as the zero scheme of a regular section 
\begin{equation*}
\sigma :P\times \Delta \rightarrow \pi ^{*}\frak{E}
\end{equation*}
where 
\begin{equation*}
\pi :P\times \Delta \rightarrow P
\end{equation*}
is the standard projection. (In what follows, we will be working in the
context of branched covers $\tilde{\Delta }\rightarrow \Delta $.) Then 
\begin{equation*}
\left. \pi ^{*}\frak{E}\right| _{X_{\Delta }}=N_{X_{\Delta }\backslash
P\times \Delta }.
\end{equation*}
Our final critical assumption is that 
\begin{equation}
\omega _{X_{\Delta }/\Delta }=\mathcal{O}_{X_{\Delta }}.  \label{ca}
\end{equation}

Let $J^{\prime }$ denote a connected Zariski-open subscheme of the Hilbert
scheme of curves on $P$ with universal curve 
\begin{equation}
\begin{array}{ccc}
J\subset J^{\prime }\times P & \overset{q}{\longrightarrow } & P \\ 
\downarrow p &  &  \\ 
J^{\prime } &  & 
\end{array}
\label{a}
\end{equation}
We assume:

i) For all $\left\{ Z\right\} \in J^{\prime }$, $Z$ is a connected smooth
subscheme of $P$ of dimension $1$.

ii) The natural map

\begin{equation}
T_{J^{\prime }}\rightarrow p_{*}N_{J\backslash J^{\prime }\times P}
\label{a123}
\end{equation}
is an isomorphism and 
\begin{equation}
R^{1}p_{*}\left( N_{J\backslash J^{\prime }\times P}\right) =0.  \label{x01}
\end{equation}

A corollary of ii) is that 
\begin{equation}
R^{1}p_{*}\left( q^{*}\frak{E}\right) .  \label{x00}
\end{equation}
Abusing notation we write 
\begin{equation*}
\begin{array}{ccc}
J\times \Delta \subset J^{\prime }\times P\times \Delta & \overset{q}{%
\longrightarrow } & P\times \Delta \\ 
\downarrow p &  &  \\ 
J^{\prime }\times \Delta &  & 
\end{array}
\end{equation*}
for the maps induced by $\left( \ref{a}\right) $ and the product structure.
Thus we have that 
\begin{equation*}
E=p_{*}q^{*}\frak{E}
\end{equation*}
is a vector bundle 
\begin{equation*}
\varepsilon :E\rightarrow J^{\prime }\times \Delta
\end{equation*}
with a distinguished section 
\begin{equation*}
s=:p_{*}q^{*}\sigma :J^{\prime }\times \Delta \rightarrow E.
\end{equation*}

\begin{lemma}
\label{BF0}The zero-scheme 
\begin{equation*}
I^{\prime }\subseteq J^{\prime }\times \Delta 
\end{equation*}
of $s$ is the intersection of the Hilbert scheme of $X_{\Delta }$ with $%
J^{\prime }\times \Delta $.
\end{lemma}

\begin{proof}
The assertion is exactly Theorem 1.5 of \cite{Kl} where $\left( \ref{x00}
\right) $ is the needed hypothesis.
\end{proof}

\subsection{The result\label{p}}

Our goal is to prove that, in the above setting, the only obstructions to
finding a continuous family of curves $I_{t}^{\prime }$ in $X_{t}$ are
cohomological. That is, suppose that 
\begin{equation*}
M_{0}\subseteq \mathcal{J}_{1}\left( X_{0}\right) =\frac{\left( H^{0}\left(
\Omega _{X_{0}}^{3}\right) +H^{1}\left( \Omega _{X_{0}}^{2}\right) \right)
^{\vee }}{H_{3}\left( X_{0};\Bbb{Z}\right) }
\end{equation*}
is the maximal abelian subvariety of $\mathcal{J}_{1}\left( X_{0}\right) $
in the annihilator of $H^{0}\left( \Omega _{X_{0}}^{3}\right) $ and the
relative intermediate Jacobian 
\begin{equation*}
\mathcal{J}_{1}\left( X_{\Delta }/\Delta \right) =\frac{\left( H^{0}\left(
\Omega _{X_{\Delta }/\Delta }^{3}\right) +H^{1}\left( \Omega _{X_{\Delta
}/\Delta }^{2}\right) \right) ^{\vee }}{H_{3}\left( X_{\Delta }/\Delta ;\Bbb{%
Z}\right) }
\end{equation*}
admits an isogeny 
\begin{equation}
\mathcal{J}_{1}\left( X_{\Delta }/\Delta \right) \rightarrow M\oplus M^{\bot
}  \label{splitset}
\end{equation}
extending the isogeny 
\begin{equation*}
\mathcal{J}_{1}\left( X_{0}\right) \rightarrow M_{0}\oplus M_{0}^{\bot }.
\end{equation*}
Then the relative Hilbert scheme $I^{\prime }/\Delta $ of $X/\Delta $ is
such that the closure of 
\begin{equation*}
\left. I^{\prime }\right| _{\Delta -\left\{ 0\right\} }
\end{equation*}
contains 
\begin{equation*}
I_{0}^{\prime }=I^{\prime }\cap \left( J^{\prime }\times \left\{ 0\right\}
\right) .
\end{equation*}
A surprising corollary is that a continuous family of rational curves on a
complete-intersection Calabi-Yau threefold $X_{0}$ whose intermediate
Jacobian has no abelian subvarieties \textit{always} deforms.

\subsection{Organization of the proof}

In the above setting, if we have the isogeny $\left( \ref{splitset}\right) $%
, we will show that, for any map 
\begin{equation*}
\tilde{\Delta }\rightarrow X^{\prime }
\end{equation*}
and $n$-th order extension $W_{n}^{\prime }$ of a small open analytic subset 
$W_{0}^{\prime }$ in the smooth locus of $\left( I_{0}^{\prime }\right)
_{red.}$ into the relative Hilbert scheme of $X_{\tilde{\Delta }_{n}}/\tilde{%
\Delta }_{n}$, all sections 
\begin{equation*}
\omega \in H^{3,0}\left( X_{\tilde{\Delta }}/\tilde{\Delta }\right)
\end{equation*}
go to zero under the composition 
\begin{equation*}
H^{3,0}\left( X_{\tilde{\Delta }}/\tilde{\Delta }\right) \overset{\nabla }{%
\longrightarrow }H^{2,1}\left( X_{\tilde{\Delta }}/\tilde{\Delta }\right) 
\overset{\left( pull-back\right) }{\longrightarrow }R^{1}p_{*}^{\prime
}\left( \Omega _{W_{n}/\tilde{\Delta }_{n}}^{2}\right) \rightarrow \Omega
_{W_{n}^{\prime }}^{1}.
\end{equation*}
We will show this by using the characterization of $I^{\prime }$ as a
gradient scheme, which is proved in \cite{C2} to prove that, for the
associated Abel-Jacobi mapping $\varphi $, 
\begin{equation*}
\varphi \left( W_{n}^{\prime }\right) \subseteq M_{n}
\end{equation*}
where $M_{n}$ is the abelian variety summand of 
\begin{equation*}
\mathcal{J}_{1}\left( X_{\tilde{\Delta }_{n}}/\tilde{\Delta }_{n}\right) .
\end{equation*}

Thus, by Theorem \ref{k}, if the $\varsigma \in R^{1}p_{*}^{\prime }\left(
N_{W_{0}\backslash W_{0}^{\prime }\times X_{0}}\right) $ is the
(curvilinear) obstruction to extending $W_{n}/W_{n^{\prime }}$ to a family
of curves in $X_{\tilde{\Delta }_{n+1}}/\tilde{\Delta }_{n+1}$, the section 
\begin{equation*}
\int\nolimits_{W_{0}/W_{0}^{\prime }}\left. \left\langle \left. \varsigma
\right| \omega _{0}\right\rangle \right| _{W_{0}}\in \Omega _{W_{0}^{\prime
}}^{1}\otimes \frac{_{\left\{ \tilde{t}^{n+1}\right\} }}{_{\left\{ \tilde{t}%
^{n+2}\right\} }}
\end{equation*}
vanishes. On the other hand, we show that the map 
\begin{eqnarray*}
R^{1}p_{*}^{\prime }\left( N_{W_{0}\backslash W_{0}^{\prime }\times
X_{0}}\right) &\rightarrow &R^{1}p_{*}^{\prime }\left( \Omega
_{W_{0}}^{2}\right) \rightarrow \Omega _{W_{0}^{\prime }}^{1} \\
\varsigma ^{\prime } &\mapsto &\int\nolimits_{W_{0}/W_{0}^{\prime
}}\left\langle \left. \varsigma ^{\prime }\right| \omega _{0}\right\rangle
\end{eqnarray*}
is surjective for any non-zero $\omega _{0}$. So our (curvilinear)
obstructions $\varsigma $ to extension all lie in a proper sub-bundle of $%
R^{1}p_{*}^{\prime }\left( N_{W_{0}\backslash W_{0}^{\prime }\times
X_{0}}\right) $ of corank equal to $\dim W_{0}^{\prime }$. This gives an
upper bound on the minimum number of equations defining $I^{\prime }$ in $%
J^{\prime }\times \Delta $ that will imply the desired result.

\section{Curvilinear obstructions}

\subsection{Kuranishi theory}

Suppose, as in $\left( \ref{nf1}\right) $, we have families of submanifolds 
\begin{equation*}
W_{n}/W_{n}^{\prime }
\end{equation*}
of 
\begin{equation*}
X_{n}/\tilde{\Delta }_{n}.
\end{equation*}
As in the proof of Theorem 12.1 of \cite{C1}, we construct a $C^{\infty }$
``transversely holomorphic'' trivialization 
\begin{equation*}
F:W_{0}^{\prime }\times X_{\tilde{\Delta }}\rightarrow W_{0}^{\prime }\times
X_{0}\times \tilde{\Delta },
\end{equation*}
which is ``adapted to $\left( \ref{nf1}\right) $,'' that is, 
\begin{equation*}
p_{n}^{-1}\left( \left\{ w\right\} \times \tilde{\Delta }_{n}\right)
\subseteq F^{-1}\left( Z_{w}\times \left\{ w\right\} \times \tilde{\Delta }%
\right)
\end{equation*}
for each $w\in $ $W_{0}^{\prime }$. Let 
\begin{equation*}
\xi =\sum\nolimits_{j>0}\xi _{j}\tilde{t}^{j}
\end{equation*}
be the Kuranishi data associated to this trivialization of the deforemation $%
X_{\tilde{\Delta }}/\tilde{\Delta }$.

We consider the associated family of complex tori 
\begin{equation*}
\mathcal{J}_{q}\left( X_{\tilde{\Delta }}/\tilde{\Delta }\right) =\frac{%
\left( H^{0}\left( \Omega _{X_{\tilde{\Delta }}/\tilde{\Delta }%
}^{2q+1}\right) +\ldots +H^{q+1}\left( \Omega _{X_{\tilde{\Delta }}/\tilde{%
\Delta }}^{q}\right) \right) ^{\vee }}{H_{2q+1}\left( X_{\tilde{\Delta }}/%
\tilde{\Delta };\Bbb{Z}\right) }.
\end{equation*}
If, as in \cite{C1}, we define 
\begin{equation*}
K^{p,q}=H^{q}\left( A_{X_{0}}^{p,*},\overline{\partial }_{X_{0}}-L_{\xi
}^{1,0}\right) ,
\end{equation*}
then in \cite{C1} it is shown that 
\begin{equation*}
H^{q}\left( \Omega _{X_{\tilde{\Delta }}/\tilde{\Delta }}^{p}\right)
=e^{\left\langle \left. \xi \right| \ \right\rangle }K^{p,q}
\end{equation*}
and so we can rewrite $\mathcal{J}_{q}\left( X_{\tilde{\Delta }}/\tilde{%
\Delta }\right) $ as 
\begin{equation*}
\mathcal{J}_{q}\left( X_{\tilde{\Delta }}/\tilde{\Delta }\right) =\frac{%
\left( e^{\left\langle \left. \xi \right| \ \right\rangle }\left(
K^{2q+1,0}+\ldots +K^{q+1,q}\right) \right) ^{\vee }}{H_{2q+1}\left( X_{0};%
\Bbb{Z}\right) }.
\end{equation*}
In this last formulation the Abel-Jacobi map 
\begin{equation*}
\varphi :W_{0}^{\prime }\times \tilde{\Delta }_{n}\rightarrow \mathcal{J}%
_{q}\left( X_{n}/\tilde{\Delta }_{n}\right)
\end{equation*}
is given by the mapping 
\begin{equation}
\left( w,t\right) \mapsto \left(
\int\nolimits_{Z_{0}}^{Z_{w}}:e^{\left\langle \left. \xi \right| \
\right\rangle }\left( K^{2q+1,0}+\ldots +K^{q+1,q}\right) \rightarrow \Bbb{C}%
\right) .  \label{goodmap}
\end{equation}

Recapping the proof of Theorem 13.1 of \cite{C1}, 
\begin{equation*}
\xi =\sum\nolimits_{j>0}\xi _{j}\tilde{t}^{j}
\end{equation*}
has by construction the property that, for $j\leq n$, 
\begin{equation}
\left. \xi _{j}\right| _{W_{0}}\in A_{W_{0}}^{0,1}\left( T_{p_{0}}\right)
\label{goodextra}
\end{equation}
where $T_{p_{0}}$ is the sheaf of tangent vectors to the fibers of the
fibration $p_{0}$. Again following \cite{C1} the action of the Gauss-Manin
connection is given by the operator 
\begin{equation*}
\left\langle \left. \frac{\partial \xi }{\partial \tilde{t}}\right| \
\right\rangle .
\end{equation*}
Thus by type $\left( \ref{goodmap}\right) $ restricted to 
\begin{equation*}
\left\langle \left. \frac{\partial \xi }{\partial \tilde{t}}\right| \
\right\rangle e^{\left\langle \left. \xi \right| \ \right\rangle }\left(
K^{2q+1,0}+\ldots +K^{q+2,q-1}\right)
\end{equation*}
is given over $\tilde{\Delta }_{n}$ by 
\begin{equation}
\int\nolimits_{Z_{0}}^{Z_{w}}:\left\langle \left. \frac{\partial \xi }{%
\partial t}\right| H^{q+2,q-1}\left( X_{0}\right) \right\rangle \rightarrow 
\Bbb{C}  \label{above}
\end{equation}
since, using $\left( \ref{goodextra}\right) $, all other summands are zero
through order $n$. The mapping $\left( \ref{above}\right) $ vanishes on $%
W_{0}^{\prime }\times \tilde{\Delta }_{n-1}$ by construction and the proof
concludes by showing that the element in 
\begin{equation*}
R^{q}p_{*}^{\prime }\left( \Omega _{W_{0}}^{p+q}\right) \otimes \frac{%
_{\left\{ \tilde{t}^{n+1}\right\} }}{_{\left\{ \tilde{t}^{n+2}\right\} }}
\end{equation*}
given by $\left( \ref{above}\right) $ is $n$ times the element \ref{BF-1}
given by capping with the obstruction class $\varsigma $ to extending the
family $W_{n}/$ $W_{n}^{\prime }$ to a family of submanifolds over $\Delta
_{n+1}$.

\subsection{Curvilinear obstructions suffice}

Since we will need to work with obstructions arising from only curvilinear
deformations, we will need the following general lemma before proceeding.
Let 
\begin{equation*}
X\rightarrow X^{\prime }
\end{equation*}
be a family of complex manifolds parametrized by a smooth space $X^{\prime }$
with distinguished point $\left\{ X_{0}\right\} $.

\begin{lemma}
\label{i}Let $Z_{0}$ be a submanifold of $X_{0}$ and let 
\begin{equation*}
\begin{tabular}{ccccc}
$Z_{0}$ & $\subseteq $ & $I$ & $\subseteq $ & $X$ \\ 
$\downarrow $ &  & $\downarrow $ &  & $\downarrow $ \\ 
$\left\{ 0\right\} $ & $\in $ & $I^{\prime }$ & $\subseteq $ & $X^{\prime }$%
\end{tabular}
\end{equation*}
be a maximal flat deformation of $Z_{0}$ in $X/X^{\prime }$. Suppose that,
for all $n>0$ and for all flat deformation 
\begin{equation*}
\begin{tabular}{ccccc}
$Z_{n}=\tilde{\Delta }_{n}\times _{I^{\prime }}I$ & $\rightarrow $ & $X_{%
\tilde{\Delta }}$ & $\rightarrow $ & $X$ \\ 
$\downarrow $ &  & $\downarrow $ &  & $\downarrow $ \\ 
$\tilde{\Delta }_{n}$ & $\rightarrow $ & $\tilde{\Delta }$ & $\rightarrow $
& $X^{\prime }$%
\end{tabular}
\end{equation*}
of $Z_{0}$, the obstructions $\upsilon $ to extending the family to a
diagram 
\begin{equation*}
\begin{tabular}{ccccc}
$Z_{n+1}$ & $\rightarrow $ & $X_{\tilde{\Delta }}$ & $\rightarrow $ & $X$ \\ 
$\downarrow $ &  & $\downarrow $ &  & $\downarrow $ \\ 
$\tilde{\Delta }_{n+1}$ & $\rightarrow $ & $\tilde{\Delta }$ & $\rightarrow $
& $X^{\prime }$%
\end{tabular}
\end{equation*}
lie in some fixed subspace 
\begin{equation*}
U\subseteq H^{1}\left( N_{Z_{0}\backslash X_{0}}\right) .
\end{equation*}
Then 
\begin{equation*}
\dim _{\left\{ Z_{0}\right\} }I^{\prime }\geq \dim X^{\prime }-\dim U.
\end{equation*}
\end{lemma}

\begin{proof}
Let $\mathcal{I}^{\prime }$ denote the ideal in the local ring of $X^{\prime
}$ at $0$ which locally defines $I^{\prime }$. The obstructions to further
extension of 
\begin{equation*}
\begin{tabular}{ccc}
$I$ & $\rightarrow $ & $X$ \\ 
$\downarrow $ &  & $\downarrow $ \\ 
$I^{\prime }$ & $\subseteq $ & $X^{\prime }$%
\end{tabular}
\end{equation*}
are given by a monomorphism 
\begin{equation*}
\upsilon :\left( \frac{\mathcal{I}^{\prime }}{\frak{m}_{0,X^{\prime }}\cdot 
\mathcal{I}^{\prime }}\right) ^{\vee }\rightarrow H^{1}\left(
N_{Z_{0}\backslash X_{0}}\right) .
\end{equation*}
(See Chapter 1 of \cite{Ko}). Let $g_{1},\ldots ,g_{s}\in \mathcal{I}%
^{\prime }$ give a basis for 
\begin{equation*}
\frac{\mathcal{I}^{\prime }}{\frak{m}_{0,X^{\prime }}\cdot \mathcal{I}%
^{\prime }}
\end{equation*}
Let $r=$ codim$_{0}\left( I^{\prime },X^{\prime }\right) \leq s$ and, using
the ``Curve Selection Lemma'' (Lemma 5.4 of \cite{FM}), choose maps 
\begin{equation*}
f_{i}:\Delta \rightarrow X^{\prime },\quad i=1,\ldots ,r,
\end{equation*}
such that 
\begin{equation*}
g_{j}\circ f_{i}=0,\quad \forall j<i
\end{equation*}
and 
\begin{equation*}
g_{i}\circ f_{i}
\end{equation*}
generates the non-zero ideal $\left\{ \tilde{t}^{n_{i}}\right\}
=f_{i}^{*}\left( \mathcal{I}^{\prime }\right) $. By construction (see, for
example, the proof of Proposition 4.7 of \cite{BF}), each element 
\begin{equation}
\upsilon \left( \left( f_{i}\right) _{*}\left( \left( \tilde{t}
^{n_{i}}\right) ^{\vee }\right) \right)  \label{j}
\end{equation}
is non-zero in the quotient vector space 
\begin{equation*}
\frac{H^{1}\left( N_{Z_{0}\backslash X_{0}}\right) }{\left\langle \upsilon
\left( \left( f_{j}\right) _{*}\left( \left( \tilde{t}^{n_{j}}\right) ^{\vee
}\right) \right) \right\rangle _{j<i}}.
\end{equation*}
So the elements $\left( \ref{j}\right) $ are linearly independent in 
\begin{equation*}
U\subseteq H^{1}\left( N_{Z_{0}\backslash X_{0}}\right) .
\end{equation*}
So $r\leq \dim U$.
\end{proof}

\section{Local analytic obstructions\label{ch3}}

\subsection{Algebraic obstruction map}

For the rest of this paper we assume that we are in the situation of \ref
{set}$.$ As in \cite{CK} we consider $E^{\vee }$ as linear functions on $E$
which we pull back to $J^{\prime }$ via the section $s=p_{*}q^{*}\sigma $.
Thus by Lemma \ref{BF0} 
\begin{equation}
s^{*}:E^{\vee }\rightarrow \mathcal{O}_{J^{\prime }\times \Delta }  \label{e}
\end{equation}
has image equal to the ideal sheaf 
\begin{equation*}
\mathcal{I}^{\prime }
\end{equation*}
of $I^{\prime }$ in $J^{\prime }\times \Delta $. we have 
\begin{equation}
\begin{tabular}{ccccc}
$E^{\vee }$ & $\overset{s^{*}}{\longrightarrow }$ & $\mathcal{I}^{\prime }$
& $\rightarrow $ & $0$ \\ 
& $\searrow $ & $\downarrow $ &  &  \\ 
&  & $\left. \Omega _{J^{\prime }\times \Delta /\Delta }^{1}\right|
_{I^{\prime }}$ &  & 
\end{tabular}
\label{-a}
\end{equation}
and so the restriction induced diagram 
\begin{equation}
\begin{tabular}{ccccc}
$E_{\Delta \times J^{\prime }}^{\vee }$ & $\overset{s_{0}^{*}}{%
\longrightarrow }$ & $\mathcal{I}^{\prime }$ & $\rightarrow $ & $0$ \\ 
& $\searrow $ & $\downarrow $ &  &  \\ 
&  & $\left. \Omega _{\Delta \times J^{\prime }}^{1}\right| _{I_{0}^{\prime
}}$ &  & 
\end{tabular}
.  \label{-ar}
\end{equation}

Let $S_{0}^{\prime }$ denote a reduced smooth quasi-projective scheme in $%
I_{0}^{\prime }$ such that 
\begin{equation*}
p_{*}N_{S_{0}\backslash S_{0}^{\prime }\times X_{0}}
\end{equation*}
is locally free. Using this diagram 
\begin{equation*}
\begin{array}{ccccccc}
E_{\left\langle S_{0}^{\prime }\right\rangle }^{\vee } & \rightarrow & 
\left. \Omega _{J^{\prime }\times \Delta }\right| _{\left\langle
S_{0}^{\prime }\right\rangle } & \rightarrow & \left. \Omega _{I^{\prime
}}\right| _{\left\langle S_{0}^{\prime }\right\rangle } & \rightarrow & 0 \\ 
\downarrow &  & \downarrow &  & \downarrow &  &  \\ 
E_{\left\langle S_{0}^{\prime }\right\rangle }^{\vee } & \rightarrow & 
\left. \Omega _{J^{\prime }}\right| _{\left\langle S_{0}^{\prime
}\right\rangle } & \rightarrow & \left. \Omega _{I_{0}^{\prime }}\right|
_{\left\langle S_{0}^{\prime }\right\rangle } & \rightarrow & 0
\end{array}
\end{equation*}
with exact rows and vertical surjections as well as $\left( \ref{-a}\right) $
and $\left( \ref{-ar}\right) $, we can choose 
\begin{equation*}
g_{1},\ldots ,g_{r^{\prime }}\in \left. \mathcal{I}^{\prime }\right|
_{\left\langle S_{0}^{\prime }\right\rangle }
\end{equation*}
such that 
\begin{equation*}
dg_{1},\ldots ,dg_{r^{\prime }}
\end{equation*}
restrict to a basis for the kernel of the morphism 
\begin{equation*}
\left. \Omega _{J^{\prime }}^{1}\right| _{\left\langle S_{0}^{\prime
}\right\rangle }\rightarrow \left. \Omega _{I_{0}^{\prime }}^{1}\right|
_{\left\langle S_{0}^{\prime }\right\rangle }.
\end{equation*}
The $g_{i}$ define an analytic space 
\begin{equation}
U_{\left\langle S_{0}^{\prime }\right\rangle }^{\prime }\subseteq \Delta
\times J^{\prime }  \label{X'}
\end{equation}
which is smooth over $\Delta $ and such that $I^{\prime }\subseteq $ $%
W_{\Delta }^{\prime }$. Now 
\begin{equation*}
\dim U_{\left\langle S_{0}^{\prime }\right\rangle }^{\prime }=\dim J^{\prime
}-r^{\prime }+1
\end{equation*}
and the dimension of the fiber of $U_{\left\langle S_{0}^{\prime
}\right\rangle }^{\prime }/\Delta $ is the rank of $p_{*}N_{S_{0}\backslash
S_{0}^{\prime }\times X_{0}}$.

Let 
\begin{equation*}
\mathcal{S}_{0}^{\prime }
\end{equation*}
be the ideal sheaf of $S_{0}^{\prime }$ (on the appropriate open dense
subset of $J^{\prime }).$ Shrinking the open dense subset if necessary, we
can assume that 
\begin{equation*}
\frac{\mathcal{I}^{\prime }}{\mathcal{S}_{0}^{\prime }\cdot \mathcal{I}%
^{\prime }}
\end{equation*}
is also locally free. Define 
\begin{equation*}
S_{0}=S_{0}^{\prime }\times _{I_{0}^{\prime }}I_{0}.
\end{equation*}
Applying $Rp_{*}$ to the exact sequence 
\begin{equation*}
0\rightarrow N_{S_{0}\backslash S_{0}^{\prime }\times X_{0}}\rightarrow
N_{S_{0}\backslash S_{0}^{\prime }\times P}\rightarrow \left.
q^{*}N_{X_{0}\backslash P}\right| _{S_{0}}\rightarrow 0
\end{equation*}
we obtain the exact sequence 
\begin{equation}
0\rightarrow p_{*}N_{S_{0}\backslash S_{0}^{\prime }\times X_{0}}\rightarrow
\left. T_{J^{\prime }}\right| _{S_{0}^{\prime }}\rightarrow E_{S_{0}^{\prime
}}\rightarrow R^{1}p_{*}N_{S_{0}\backslash S_{0}^{\prime }\times
X_{0}}\rightarrow 0.  \label{-a'}
\end{equation}

Now from $\left( \ref{-a}\right) $ we have exact 
\begin{equation}
E_{S_{0}^{\prime }}^{\vee }\rightarrow \frac{\mathcal{I}^{\prime }}{\mathcal{%
S}_{0}^{\prime }\cdot \mathcal{I}^{\prime }}\rightarrow 0  \label{-a"}
\end{equation}
from which we achieve an injection 
\begin{equation*}
\left( \frac{\mathcal{I}^{\prime }}{\mathcal{S}_{0}^{\prime }\cdot \mathcal{I%
}^{\prime }}\right) ^{\vee }\rightarrow E_{S_{0}^{\prime }}
\end{equation*}
and so, composing with the surjection 
\begin{equation*}
E_{S_{0}^{\prime }}\rightarrow R^{1}p_{*}N_{S_{0}\backslash S_{0}^{\prime
}\times X_{0}}
\end{equation*}
in $\left( \ref{-a'}\right) $ we construct a map 
\begin{equation*}
\alpha :\left( \frac{\mathcal{I}^{\prime }}{\mathcal{S}_{0}^{\prime }\cdot 
\mathcal{I}^{\prime }}\right) ^{\vee }\rightarrow
R^{1}p_{*}N_{S_{0}\backslash S_{0}^{\prime }\times X_{0}}.
\end{equation*}

\subsection{Anaytic obstruction map}

On the other hand, by general obstruction theory we have an injective map 
\begin{equation}
\beta :\left( \frac{\mathcal{I}^{\prime }}{\mathcal{S}_{0}^{\prime }\cdot 
\mathcal{I}^{\prime }}\right) ^{\vee }\rightarrow
R^{1}p_{*}N_{S_{0}\backslash S_{0}^{\prime }\times X_{0}}  \label{-b}
\end{equation}
constructed via Kuranishi theory. (See \S 13 of \cite{C1}; also, for
example, Chapter 1 of \cite{Ko}). As in \cite{C2} we start with a local
versal deformation $X/X^{\prime }$ of $X_{0}$ a point $\left\{ Y_{0}\right\}
\in S_{0}^{\prime }$ and let$\ U^{\prime }$ be smooth over $X^{\prime }$ of
minimal dimension containing the local (relative) Hilbert scheme $Y^{\prime
} $. Then 
\begin{equation*}
Y_{\Delta }^{\prime }=Y^{\prime }\times _{X^{\prime }}\Delta \subseteq
I^{\prime }.
\end{equation*}
Furthermore we can assume that 
\begin{equation*}
U_{\Delta }^{\prime }=U^{\prime }\times _{X^{\prime }}\Delta
\end{equation*}
is an open neighborhood of $\left( \left\{ Y_{0}\right\} ,\left\{
X_{0}\right\} \right) $ in the smooth space $U_{\left\langle S_{0}^{\prime
}\right\rangle }^{\prime }$ constructed in $\left( \ref{X'}\right) $.

In \cite{C2} one constructs a $C^{\infty }$ ``transversely holomorphic''
trivialization 
\begin{equation*}
F:X\times _{X^{\prime }}U^{\prime }\rightarrow X_{0}\times U^{\prime },
\end{equation*}
which is ``adapted to the universal curve $Y/Y^{\prime }$ at $\left\{
Y_{0}\right\} $,'' that is, 
\begin{equation*}
Y\subseteq F^{-1}\left( Y_{0}\times U^{\prime }\right) .
\end{equation*}

For each analytic map 
\begin{equation*}
\vartheta :\tilde{\Delta }\rightarrow U_{\Delta }^{\prime }
\end{equation*}
sending $0$ to a point $s_{0}^{\prime }\in Y_{\Delta }\cap S_{0}^{\prime }$
and having the property that 
\begin{equation*}
\vartheta ^{*}\left( \mathcal{I}^{\prime }\right) =\left\{ \tilde{t}%
^{n+1}\right\} ,
\end{equation*}
the Kuranishi data associated to the trivialization $F$ produces the
obstruction element 
\begin{equation*}
\varsigma \in H^{1}\left( N_{S_{S_{0}^{\prime }}\backslash X_{0}}\right)
\end{equation*}
in $\left( \ref{genobst}\right) $ as the image of 
\begin{equation*}
\left( \frac{\tilde{t}^{n+1}}{\tilde{t}^{n+2}}\right) ^{\vee }\subseteq
\left( \frac{\mathcal{I}^{\prime }}{\mathcal{S}_{0}^{\prime }\cdot \mathcal{I%
}^{\prime }}\right) ^{\vee }
\end{equation*}
under $\beta $.

We wish to show that under certain base extensions, $\alpha $ and $\beta $
give the same obstruction map. But first we need to make precise the
pull-backs we will need for this assertion. Let 
\begin{equation*}
V=Y_{\Delta }^{\prime }\cap S_{0}^{\prime }
\end{equation*}
as above, let $D$ an auxiliary polydisk, and suppose we have an analytic
isomorphism 
\begin{equation*}
\psi :V\times D\rightarrow U_{\Delta }^{\prime }
\end{equation*}
such that 
\begin{equation*}
\left. \psi \right| _{V\times \left\{ 0\right\} }=identity_{V}.
\end{equation*}
By flat stratification we can, after shrinking $V$ if necessary, arrange
that the family of ideals 
\begin{equation*}
\left. \psi ^{*}\left( \mathcal{I}^{\prime }\right) \right| _{\left\{
v\right\} \times D}
\end{equation*}
is flat over $V$. We then use the ``Curve Selection Lemma'' (Lemma 5.4 of 
\cite{FM}) as in the proof of Lemma \ref{i} on each $\left\{ v\right\}
\times D$, but this time with analytic parameter $v$. Letting $\tilde{t}$
denote the analytic parameter of a complex disk $\tilde{\Delta }$ , we can,
after possibly again shrinking $V$, choose analytic maps 
\begin{equation}
f_{i}:V\times \tilde{\Delta }\rightarrow V\times D,\quad i=1,\ldots ,r,
\label{cursel}
\end{equation}
such that 
\begin{equation*}
g_{j}\circ f_{i}=0,\quad \forall j<i
\end{equation*}
and 
\begin{equation*}
g_{i}\circ f_{i}
\end{equation*}
generates the non-zero ideal $\left\{ \tilde{t}^{n_{i}}\right\}
=f_{i}^{*}\left( \mathcal{I}^{\prime }\right) $. By construction (see again,
for example, the proof of Proposition 4.7 of \cite{BF}), each element 
\begin{equation*}
\beta _{*}\left( \psi _{*}\left( \left( f_{i}\right) _{*}\left( \left( 
\tilde{t}^{n_{i}}\right) ^{\vee }\right) \right) \right)
\end{equation*}
gives a never-zero section of the quotient vector bundle 
\begin{equation*}
\frac{R^{1}p_{*}\left( N_{Y_{V}\backslash V\times X_{0}}\right) }{%
\left\langle \beta _{*}\left( \psi _{*}\left( \left( f_{j}\right) _{*}\left(
\left( \tilde{t}^{n_{j}}\right) ^{\vee }\right) \right) \right)
\right\rangle _{j<i}}.
\end{equation*}

\subsection{Equality of the two obstruction maps}

We claim that, for each $i$, 
\begin{equation}
\alpha _{*}\left( \psi _{*}\left( \left( f_{i}\right) _{*}\left( \left( 
\tilde{t}^{n_{i}}\right) ^{\vee }\right) \right) \right) =\beta _{*}\left(
\psi _{*}\left( \left( f_{i}\right) _{*}\left( \left( \tilde{t}%
^{n_{i}}\right) ^{\vee }\right) \right) \right)  \label{comclaim}
\end{equation}
To prove the claim, it suffices to prove the following lemma.

\begin{lemma}
\label{key}Let $\left\{ Z_{0}\right\} \in S_{0}^{\prime }$ and let 
\begin{eqnarray*}
\gamma  &:&\tilde{\Delta }\rightarrow J^{\prime }\times \Delta  \\
\tilde{t} &\mapsto &\left( j^{\prime }\left( \tilde{t}\right) ,t\left( 
\tilde{t}\right) \right) 
\end{eqnarray*}
be an analytic map of the $\tilde{t}$-disk $\tilde{\Delta }$ such that 
\begin{equation*}
\gamma \left( 0\right) =\left\{ Z_{0}\right\} .
\end{equation*}
Let 
\begin{equation*}
\left\{ \tilde{t}^{m+1}\right\} =\gamma ^{*}\left( \mathcal{I}^{\prime
}\right) .
\end{equation*}
Then 
\begin{equation*}
\left( \alpha \circ \gamma _{*}\right) \left( \frac{\left\{ \tilde{t}%
^{m+1}\right\} }{\left\{ \tilde{t}^{m+2}\right\} }\right) =\left( \beta
\circ \gamma _{*}\right) \left( \frac{\left\{ \tilde{t}^{m+1}\right\} }{%
\left\{ \tilde{t}^{m+2}\right\} }\right) .
\end{equation*}
\end{lemma}

\begin{proof}
We consider the pull-back family of threefolds 
\begin{equation*}
\begin{array}{ccc}
\tilde{X}_{\tilde{\Delta }} & \rightarrow & X_{\Delta } \\ 
\downarrow &  & \downarrow \\ 
\tilde{\Delta } & \overset{t}{\longrightarrow } & \Delta
\end{array}
\end{equation*}
given by the zero-scheme of the section $\tilde{\sigma}$ of the pull-back
bundle 
\begin{equation*}
\left( id.,t\left( \tilde{t}\right) \right) ^{*}\pi ^{*}\frak{E}
\end{equation*}
on 
\begin{equation*}
P\times \tilde{\Delta }.
\end{equation*}
We have a family of curves 
\begin{equation}
\tilde{I}_{n}/\tilde{\Delta }_{n}\subseteq \tilde{X}_{n}/\tilde{\Delta }_{n}
\label{critfam}
\end{equation}
via pull-back of the family $\left( J/J^{\prime }\right) \times \Delta $
under $\gamma $. Since 
\begin{equation*}
R^{1}p_{*}N_{J\backslash J^{\prime }\times P}=0
\end{equation*}
the family $\left( \ref{critfam}\right) $ extends to a family of curves 
\begin{equation*}
\begin{array}{ccc}
\tilde{J}\subseteq P\times \tilde{\Delta } & \overset{\tilde{q}}{%
\longrightarrow } & P \\ 
\downarrow ^{\tilde{p}} &  &  \\ 
\tilde{\Delta } &  & 
\end{array}
.
\end{equation*}
Then by Lemma \ref{BF0} the zero-scheme of the section 
\begin{equation*}
\tilde{s}=\tilde{p}_{*}\tilde{q}^{*}\tilde{\sigma}
\end{equation*}
of 
\begin{equation*}
\tilde{E}=\tilde{p}_{*}\tilde{q}^{*}\tilde{\sigma}
\end{equation*}
is exactly 
\begin{equation*}
\tilde{\Delta }_{n}.
\end{equation*}

Now construct a (transversely holomorphic) $C^{\infty }$-trivialization 
\begin{equation*}
F:P\times \tilde{\Delta }\rightarrow P\times \tilde{\Delta }
\end{equation*}
over $\tilde{\Delta }$ such that

i) 
\begin{equation*}
F\left( X\right) =X_{0}\times \tilde{\Delta },
\end{equation*}

ii) 
\begin{equation*}
F\left( \tilde{J}_{n}\right) \subseteq Z_{0}\times \tilde{\Delta }.
\end{equation*}
Then as in \S 11 of \cite{C1} the Kuranishi data 
\begin{equation*}
\tilde{\xi}=\sum\nolimits_{j>0}\tilde{\xi}_{j}\tilde{t}^{j}
\end{equation*}
associated to the trivialization $F$ has the property that 
\begin{equation*}
\left\{ \left. \tilde{\xi}_{n+1}\right| _{Z_{0}}\right\} \in H^{1}\left(
N_{Z_{0}\backslash X_{0}}\right)
\end{equation*}
is the obstruction class to extending the family $\left( \ref{critfam}
\right) $ to a family of curves in $\tilde{X}_{n+1}/\tilde{\Delta }_{n+1}$.
But by the construction of the trivialization $F$ above, this class is given
by the obstruction to extending the map 
\begin{equation*}
\tilde{J}_{n}/\tilde{\Delta }_{n}\rightarrow \tilde{X}/\tilde{\Delta }
\end{equation*}
to a map 
\begin{equation*}
\tilde{J}_{n+1}/\tilde{\Delta }_{n+1}\rightarrow \tilde{X}/\tilde{\Delta }.
\end{equation*}
But this last obstruction class, which lives in 
\begin{equation*}
H^{1}\left( N_{\Gamma _{Z_{0}}\backslash Z_{0}\times X_{0}}\right)
=H^{1}\left( \left. \tilde{q}^{*}T_{\tilde{X}_{0}}\right| _{Z_{0}}\right) ,
\end{equation*}
is the image of 
\begin{equation*}
\tilde{s}\in \frac{\left\{ \tilde{t}^{n+1}\right\} \tilde{E}}{\left\{ \tilde{%
t}^{n+2}\right\} \tilde{E}}=\tilde{E}_{\left\{ Z_{0}\right\} }
\end{equation*}
under the map 
\begin{equation*}
\tilde{E}_{\left\{ Z_{0}\right\} }\rightarrow H^{1}\left( \left.
T_{X_{0}}\right| _{Z_{0}}\right)
\end{equation*}
induced by the exact sequence 
\begin{equation*}
0\rightarrow \left. \tilde{q}^{*}T_{X_{0}}\right| _{Z_{0}}\rightarrow \left. 
\tilde{q}^{*}T_{P}\right| _{Z_{0}}\rightarrow \left. \frak{E}\right|
_{Z_{0}}\rightarrow 0.
\end{equation*}
\end{proof}

\subsection{Obstructions as K\"{a}hler differentials}

\begin{lemma}
The map 
\begin{equation}
R^{1}p_{*}N_{Y_{0}\backslash Y_{0}^{\prime }\times X_{0}}\rightarrow \Omega %
_{Y_{0}^{\prime }}^{1}  \label{smap}
\end{equation}
given in $\left( \ref{BF-1}\right) $ in Theorem \ref{k} is surjective.
\end{lemma}

\begin{proof}
Let $\Gamma ^{\prime }$ denote the diagonal of $S_{0}^{\prime }\times
S_{0}^{\prime }$ and let 
\begin{equation*}
\Gamma
\end{equation*}
be the inverse image of $\Gamma ^{\prime }$ under 
\begin{equation*}
S_{0}^{\prime }\times S_{0}\rightarrow S_{0}^{\prime }\times S_{0}^{\prime }.
\end{equation*}
Then 
\begin{eqnarray*}
\Omega _{S_{0}^{\prime }}^{1} &=&N_{\Gamma ^{\prime }\backslash
S_{0}^{\prime }\times S_{0}^{\prime }}^{*} \\
&=&R^{1}p_{*}\left( N_{\Gamma \backslash S_{0}^{\prime }\times
S_{0}}^{*}\otimes \omega _{S_{0}/S_{0}^{\prime }}\right)
\end{eqnarray*}
this latter isomorphism being often called ``integration over the fiber.''
Now consider the map 
\begin{equation*}
N_{S_{0}\backslash S_{0}^{\prime }\times X_{0}}^{*}\rightarrow
N_{S_{0}\backslash S_{0}^{\prime }\times S_{0}}^{*}
\end{equation*}
induced by 
\begin{equation*}
S_{0}^{\prime }\times S_{0}\rightarrow S_{0}^{\prime }\times X_{0}.
\end{equation*}
The induced map 
\begin{equation}
R^{1}p_{*}\left( N_{S_{0}\backslash S_{0}^{\prime }\times X_{0}}^{*}\otimes
\omega _{S_{0}/S_{0}^{\prime }}\right) \rightarrow R^{1}p_{*}\left(
N_{\Gamma \backslash S_{0}^{\prime }\times S_{0}}^{*}\otimes \omega
_{S_{0}/S_{0}^{\prime }}\right)  \label{1"BF}
\end{equation}
is surjective because its dual via Verdier duality 
\begin{equation*}
T_{S_{0}^{\prime }}\rightarrow p_{*}N_{S_{0}\backslash S_{0}^{\prime }\times
X_{0}}
\end{equation*}
is injective. Composing $\left( \ref{1"BF}\right) $ with integration over
the fiber we obtain a surjection 
\begin{equation}
R^{1}p_{*}\left( N_{S_{0}\backslash S_{0}^{\prime }\times X_{0}}^{*}\otimes
\omega _{S_{0}/S_{0}^{\prime }}\right) \rightarrow \Omega _{S_{0}^{\prime
}}^{1}.  \label{1'BF}
\end{equation}

Finally recall that $\mathcal{S}_{0}$ denotes the ideal of $S_{0}\subseteq
S_{0}^{\prime }\times X_{0}$. Since $\omega _{X_{0}}$ is trivial and $X_{0}$
is a threefold we have 
\begin{equation*}
\bigwedge\nolimits^{2}\left( \mathcal{S}_{0}/\mathcal{S}_{0}^{2}\right)
\otimes \omega _{S_{0}/S_{0}^{\prime }}=q^{*}\omega _{X_{0}}=\mathcal{O}%
_{S_{0}}.
\end{equation*}
Thus 
\begin{equation*}
N_{S_{0}\backslash S_{0}^{\prime }\times X_{0}}^{*}\otimes \omega
_{S_{0}/S_{0}^{\prime }}=N_{S_{0}\backslash S_{0}^{\prime }\times X_{0}}
\end{equation*}
or, said otherwise, 
\begin{eqnarray}
N_{S_{0}\backslash S_{0}^{\prime }\times X_{0}} &\cong &N_{S_{0}\backslash
S_{0}^{\prime }\times X_{0}}^{*}\otimes \omega _{S_{0}/S_{0}^{\prime }}
\label{2BF} \\
\nu &\mapsto &\left\langle \left. \nu \right| \alpha \right\rangle  \notag
\end{eqnarray}
where $\alpha $ is a fixed generator of $H^{0}\left( \omega _{X_{0}}\right) $%
. Composing with $\left( \ref{1'BF}\right) $ we have the desired surjection 
\begin{eqnarray}
R^{1}p_{*}N_{S_{0}\backslash S_{0}^{\prime }\times X_{0}} &\rightarrow
&R^{1}p_{*}\left( N_{S_{0}\backslash S_{0}^{\prime }\times X_{0}}^{*}\otimes
\omega _{S_{0}/S_{0}^{\prime }}\right) \rightarrow \Omega _{S_{0}^{\prime
}}^{1}.  \label{2'BF} \\
\nu &\mapsto &\left\langle \left. \nu \right| q^{*}\alpha \right\rangle 
\notag
\end{eqnarray}
\end{proof}

\subsection{The crucial estimate}

Now the restriction of the map $\left( \ref{2'BF}\right) $ to $Y_{0}^{\prime
}$ is the map $\left( \ref{BF-1}\right) $ in Theorem \ref{k}. So, in a
situation in which the hypotheses of Theorem \ref{k} are satisfied, the
obstructions $\left( \ref{comclaim}\right) $ lie in the kernel of the
surjection 
\begin{equation*}
R^{1}p_{*}N_{S_{0}\backslash S_{0}^{\prime }\times X_{0}}\rightarrow \Omega
_{S_{0}^{\prime }}^{1}.
\end{equation*}
Then applying Lemma \ref{i} we can conclude 
\begin{equation}
\dim _{\left\{ Z_{0}\right\} }I^{\prime }\geq \dim U_{\left\langle
S_{0}^{\prime }\right\rangle }^{\prime }-\left( rank\left(
R^{1}p_{*}N_{S_{0}\backslash S_{0}^{\prime }\times X_{0}}\right) -\dim
S_{0}^{\prime }\right) .  \label{BBF}
\end{equation}
We now come to the main result of this paper, which treats a situation in
which the hypotheses of Theorem \ref{k} are indeed satisfied.

\section{The Hilbert scheme as gradient variety}

Suppose now that 
\begin{equation*}
M_{0}\subseteq \mathcal{J}_{1}\left( X_{0}\right) =\frac{\left( H^{0}\left(
\Omega _{X_{0}}^{3}\right) +H^{1}\left( \Omega _{X_{0}}^{2}\right) \right)
^{\vee }}{H_{3}\left( X_{0};\Bbb{Z}\right) }
\end{equation*}
is the maximal abelian subvariety of $\mathcal{J}_{1}\left( X_{0}\right) $
in the annihilator of $H^{0}\left( \Omega _{X_{0}}^{3}\right) $ and the
relative intermediate Jacobian 
\begin{equation*}
\mathcal{J}_{1}\left( X_{\Delta }/\Delta \right) =\frac{\left( H^{0}\left(
\Omega _{X_{\Delta }/\Delta }^{3}\right) +H^{1}\left( \Omega _{X_{\Delta
}/\Delta }^{2}\right) \right) ^{\vee }}{H_{3}\left( X_{\Delta }/\Delta ;\Bbb{%
\ Z}\right) }
\end{equation*}
admits an isogeny 
\begin{equation}
\mathcal{J}_{1}\left( X_{\Delta }/\Delta \right) \rightarrow M\oplus
_{\Delta }M^{\bot }  \label{isogeny'}
\end{equation}
extending the isogeny 
\begin{equation*}
\mathcal{J}_{1}\left( X_{0}\right) \rightarrow M_{0}\oplus M_{0}^{\bot }.
\end{equation*}
Let 
\begin{equation*}
\overline{I^{\prime }}
\end{equation*}
denote the Zariski closure of $I^{\prime }$ in $J^{\prime }\times \Delta $.
The composition of the Abel-Jacobi mapping 
\begin{equation*}
\overline{I^{\prime }}\rightarrow \mathcal{J}_{1}\left( X_{\Delta }/\Delta
\right)
\end{equation*}
with projection gives a map 
\begin{equation*}
\overline{I^{\prime }}\rightarrow M^{\bot }
\end{equation*}
whose image is proper over $\Delta $. But ,by the maximality of $M$, the
image must therefore be finite over $\Delta $.

It is at this point that the gradient variety construction in \cite{C2}
enters. There we consider $X^{\prime }$, a versal local deformation of the $%
K $-trivial threefold $X_{0}$, $U^{\prime }$ smooth over $X^{\prime }$ of
minimal dimension containing the local (relative) Hilbert scheme $Y^{\prime
} $, and 
\begin{eqnarray*}
\tilde{X}^{\prime } &=&\left\{ \left( x^{\prime },\omega _{x^{\prime
}}\right) :x^{\prime }\in X^{\prime },\ \omega _{x^{\prime }}\in \left(
H^{0}\left( \Omega _{X/X^{\prime }}^{3}\right) -\left\{ 0\right\} \right)
\right\} \\
\tilde{U}^{\prime } &=&U^{\prime }\times _{X^{\prime }}\tilde{X}^{\prime }%
\overset{\tilde{\pi }}{\longrightarrow }\tilde{X}^{\prime }.
\end{eqnarray*}
We choose a holomorphic product structure 
\begin{eqnarray*}
U^{\prime } &=&U_{0}^{\prime }\times X^{\prime } \\
\tilde{U}^{\prime } &=&U_{0}^{\prime }\times \tilde{X}^{\prime }
\end{eqnarray*}
and use it to define a map 
\begin{equation*}
\kappa :U_{\Delta }^{\prime }=U^{\prime }\times _{X^{\prime }}\Delta
\rightarrow \tilde{U}^{\prime }
\end{equation*}
which lifts the natural map 
\begin{equation*}
U_{\Delta }^{\prime }\rightarrow U^{\prime }.
\end{equation*}
Theorem 6.1 of \cite{C2} gives a holomorphic function $\Phi $ on $\tilde{U}%
^{\prime }$ such that, with respect to the exact sequence 
\begin{equation*}
0\rightarrow \tilde{\pi }^{*}\Omega _{\tilde{X}^{\prime }}^{1}\rightarrow
\Omega _{\tilde{U}^{\prime }}^{1}\rightarrow \Omega _{\tilde{U}^{\prime }/%
\tilde{X}^{\prime }}^{1}\rightarrow 0,
\end{equation*}
we have:

\begin{description}
\item[Property 1]  The relative Hilbert scheme $\tilde{Y}^{\prime }$,
considered as an analytic subscheme of $\tilde{U}^{\prime }$ is the
zero-scheme of the section 
\begin{equation*}
d_{\tilde{U}^{\prime }/\tilde{X}^{\prime }}\Phi
\end{equation*}
of 
\begin{equation*}
\Omega _{\tilde{U}^{\prime }/\tilde{X}^{\prime }}^{1}.
\end{equation*}

\item[Property 2]  Under a natural isomorphism 
\begin{eqnarray*}
F^{2}H^{3}\left( \tilde{X}/\tilde{X}^{\prime }\right) &\cong &T_{\tilde{X}%
^{\prime }} \\
\Omega _{\tilde{X}^{\prime }}^{1} &\cong &\left( F^{2}H^{3}\left( \tilde{X}/%
\tilde{X}^{\prime }\right) \right) ^{\vee }
\end{eqnarray*}
given by Donagi-Markman, the section 
\begin{equation*}
\left. d\Phi \right| _{\tilde{Y}^{\prime }}
\end{equation*}
of 
\begin{equation*}
\tilde{\pi }^{*}\Omega _{\tilde{X}^{\prime }}^{1}
\end{equation*}
is the normal function 
\begin{equation}
\int\nolimits_{r\left( L\times _{X^{\prime }}Y^{\prime }\right) /Y^{\prime
}}^{Y/Y^{\prime }}:\tilde{Y}^{\prime }\rightarrow \left( F^{2}H^{3}\left( 
\tilde{X}/\tilde{X}^{\prime }\right) \right) ^{\vee }.  \label{nlfn}
\end{equation}
\end{description}

Using the product structure on $\tilde{U}^{\prime }$ therefore, the
holomorphic section $\tilde{\pi }^{*}d_{\tilde{X}^{\prime }}\Phi $ of the
bundle $\left( F^{2}H^{3}\left( X\times _{X^{\prime }}U^{\prime }/U^{\prime
}\right) \right) ^{\vee }$ extend the Abel-Jacobi mapping $\left( \ref{nlfn}
\right) $. Thus the pullback $\kappa ^{*}\left( \tilde{\pi }^{*}d_{\tilde{X}%
^{\prime }}\Phi \right) $ gives a holomorphic mapping 
\begin{equation}
U_{\Delta }^{\prime }\rightarrow \left( F^{2}H^{3}\left( X_{U_{\Delta
}^{\prime }}/U_{\Delta }^{\prime }\right) \right) ^{\vee }  \label{AJ'}
\end{equation}
extending the Abel-Jacobi mapping on $Y_{\Delta }^{\prime }$.

As above let 
\begin{equation*}
V=Y_{\Delta }^{\prime }\cap S_{0}^{\prime }.
\end{equation*}
Suppose, for a mappong of disks 
\begin{equation*}
\gamma :\tilde{\Delta }\rightarrow \Delta ,
\end{equation*}
such that $\gamma \left( 0\right) =0$, we have a mapping 
\begin{equation*}
\delta :V\times \tilde{\Delta }\rightarrow U_{\Delta }^{\prime }
\end{equation*}
such that 
\begin{equation*}
Y_{\tilde{\Delta }_{n}}/Y_{\tilde{\Delta }_{n}}^{\prime }
\end{equation*}
is an extension to a family of curves in $X_{\tilde{\Delta }_{n}}$. Then
there is a holomorphic map 
\begin{equation*}
\delta :V\times \tilde{\Delta }\rightarrow U_{\Delta }^{\prime }
\end{equation*}
such that

1) 
\begin{equation*}
\delta \left( s_{0}^{\prime },0\right) =\left( \left\{ S_{s_{0}^{\prime
}}\right\} ,\left\{ X_{0}\right\} \right)
\end{equation*}

2) for all $\tilde{t}\in \tilde{\Delta }$, $\delta $ immerses $V\times
\left\{ \tilde{t}\right\} $ into the fiber of $X^{\prime }$ over $\left\{
\gamma \left( t\right) \right\} $,

3) 
\begin{equation*}
\delta \left( V\times \tilde{\Delta }_{n}\right) =Y_{\tilde{\Delta }%
_{n}}^{\prime }.
\end{equation*}

Pulling back by $\delta $, we have from $\left( \ref{AJ'}\right) $ an
extension 
\begin{equation}
V\times \tilde{\Delta }\rightarrow \left( F^{2}H^{3}\left( X_{\tilde{\Delta }%
}/\tilde{\Delta }\right) \right) ^{\vee }  \label{ext}
\end{equation}
of the Abel-Jacobi map on 
\begin{equation*}
V\times \tilde{\Delta }_{n}.
\end{equation*}

Now by hypothesis we have an isogeny 
\begin{equation}
\mathcal{J}_{1}\left( X_{\tilde{\Delta }}/\tilde{\Delta }\right) \rightarrow 
\tilde{M}\oplus _{\tilde{\Delta }}\tilde{M}^{\bot }.  \label{bigisog}
\end{equation}
so 
\begin{equation*}
\left( F^{2}H^{3}\left( X_{\tilde{\Delta }}/\tilde{\Delta }\right) \right)
^{\vee }\cong \tilde{\Delta }\times \Bbb{C}^{r_{1}}\times \Bbb{C}^{r_{2}}
\end{equation*}
where the vector spaces $\Bbb{C}^{r_{1}}$ and $\Bbb{C}^{r_{2}}$ correspond
to the summands $\tilde{M}$ and $\tilde{M}^{\bot }$ respectively. So there
is a projection of $\left( \ref{ext}\right) $ onto the second factor given
by $r_{2}$ holomorphic functions 
\begin{equation*}
g_{i}\left( s_{0}^{\prime },t\right) =\sum\nolimits_{k=0}^{\infty
}g_{ik}\left( s_{0}^{\prime }\right) \tilde{t}^{k}.
\end{equation*}
where each $g_{ik}$ is a holomorphic function on $V$. Suppose for some $i$
that $g_{ik}$ is not a constant function for some $k\leq n$. Then the
projection of the Abel-Jacobi map 
\begin{equation*}
\varphi :V\times \tilde{\Delta }_{n}\rightarrow \mathcal{J}_{1}\left( X_{%
\tilde{\Delta }_{n}}/\tilde{\Delta }_{n}\right)
\end{equation*}
into the factor $\tilde{M}_{n}^{\bot }:=\tilde{\Delta }_{n}\times _{\tilde{
\Delta }}\tilde{M}^{\bot }$ cannot lie inside an finite analytic subscheme
of $\tilde{M}_{m}^{\bot }$. However 
\begin{equation*}
Y_{\tilde{\Delta }_{n}}^{\prime }\subseteq I^{\prime }\times _{\Delta
}\Delta _{n}^{\prime }
\end{equation*}
and, as we have seen above, the image of the composition 
\begin{equation*}
I^{\prime }\times _{\Delta }\Delta _{n}^{\prime }\rightarrow \mathcal{J}%
_{1}\left( X_{\tilde{\Delta }_{n}}/\tilde{\Delta }\right) \rightarrow \tilde{%
M}_{n}^{\bot }
\end{equation*}
must be finite. Thus the functions $g_{ik}\left( s_{0}^{\prime }\right) $
are constant for all $k\leq n$. Adjusting the isogeny $\left( \ref{bigisog}%
\right) $ by a translation by these constants, we have 
\begin{equation*}
\varphi \left( V\times \tilde{\Delta }_{n}\right) \subseteq \tilde{M}.
\end{equation*}
So by Corollary \ref{goodcor}, the pairing 
\begin{equation*}
R^{1}p_{*}^{\prime }\left( \Omega _{V}^{3}\right) \rightarrow \Omega _{V}^{1}
\end{equation*}
given by $\left( \ref{BF-2}\right) $ is zero.

\section{The main theorem}

Continuing with the situation of \ref{set}, for 
\begin{equation*}
I_{0}^{\prime }:=\left( I^{\prime }\cap \left( \left\{ 0\right\} \times
J^{\prime }\right) \right) ,
\end{equation*}
let $S_{0}^{\prime }$ be a component of 
\begin{equation*}
\left( I_{0}^{\prime }\right) _{red}
\end{equation*}
and let $\left\langle S_{0}^{\prime }\right\rangle $ denote the generic
point of $S_{0}^{\prime }$. We let $\left\langle S_{0}\right\rangle $ denote 
$p^{-1}\left( \left\langle S_{0}^{\prime }\right\rangle \right) $ where 
\begin{equation*}
p:S_{0}\rightarrow S_{0}^{\prime }
\end{equation*}
is the universal curve$.$ Then 
\begin{equation*}
\left. R^{1}p_{*}\left( Hom_{I}\left( \mathcal{I}/\mathcal{I}^{2},\mathcal{O}%
_{I}\right) \right) \right| _{\left\langle S_{0}^{\prime }\right\rangle }
\end{equation*}
is a vector space over $\left\langle S_{0}^{\prime }\right\rangle $ of
dimension 
\begin{equation*}
h^{1}\left( N_{Z_{0}\backslash X_{0}}\right)
\end{equation*}
where, referring to \S \ref{ch3}$,$ $\left\{ Z_{0}\right\} \in Y_{0}^{\prime
}\subset S_{0}$. We consider the Abel-Jacobi map 
\begin{equation}
\varphi _{0}:S_{0}^{\prime }\rightarrow A_{0}\subseteq \mathcal{J}_{1}\left(
X_{0}\right) =\frac{\left( H^{0}\left( \Omega _{X_{0}}^{3}\right)
+H^{1}\left( \Omega _{X_{0}}^{2}\right) \right) ^{\vee }}{H_{3}\left( X_{0};%
\Bbb{Z}\right) }  \label{f}
\end{equation}
where $A_{0}$ is an abelian variety spanned by $\varphi _{0}\left(
S_{0}^{\prime }\right) $. This map is determined by the choice of basepoint $%
\left\{ Z_{0}\right\} \in S_{0}^{\prime }$. We write 
\begin{equation*}
A_{0}\subseteq M_{0}\subseteq \mathcal{J}_{1}\left( X_{0}\right)
\end{equation*}
where $M_{0}$ denotes the maximal abelian subvariety (lying in $H^{0}\left(
\Omega _{X_{0}}^{3}\right) ^{\bot }$).

\begin{theorem}
\label{g} Given $\left( \ref{f}\right) $, suppose that 
\begin{equation*}
M/\Delta \subseteq \mathcal{J}_{1}\left( X/\Delta \right) =\frac{\left(
H^{0}\left( \Omega _{X/\Delta }^{3}\right) +H^{1}\left( \Omega _{X/\Delta
}^{2}\right) \right) ^{\vee }}{H_{3}\left( X/\Delta ;\Bbb{Z}\right) }
\end{equation*}
is an abelian subvariety over $\Delta $ extending $M_{0}$. Then (possibly
after base extension) there is a family of curves on generic $X_{t}$
generating an abelian subvariety $A_{t}^{\prime }\subseteq \mathcal{J}%
_{1}\left( X_{t}\right) $ via the Abel-Jacobi mapping such that $%
A_{0}\subseteq A_{0}^{\prime }\subseteq M_{0}$.
\end{theorem}

\begin{proof}
We first claim that the inequality $\left( \ref{BBF}\right) $ holds. To
prove this, we must check that the hypotheses of Theorem \ref{k} hold for
one-parameter local \textit{analytic} families extending $%
Y_{0}/Y_{0}^{\prime }$. But, since we are in a projective algebraic
situation, these analytic families must come from analytic localizations as
in \S \ref{ch3} of maximal algebraic families. Using a small disk $\tilde{%
\Delta }$ with parameter $\tilde{t}$ as before, these algebraic families $%
T^{\prime }/\tilde{\Delta }$ are of the following type. We consider a
morphism 
\begin{equation}
\begin{array}{ccc}
T^{\prime } & \overset{\rho }{\longrightarrow } & \Delta \times J^{\prime }
\\ 
\downarrow &  & \downarrow \\ 
\tilde{\Delta } & \rightarrow & \Delta
\end{array}
\label{y}
\end{equation}
such that $T^{\prime }/\tilde{\Delta }$ is proper and flat and, if $%
T_{0}^{\prime }$ denotes the fiber over $0\in \tilde{\Delta }$, $\dim _{\Bbb{%
C}}T_{0}^{\prime }=\dim _{\Bbb{C}}S_{0}^{\prime }$ and for some generic
point $\left\langle T_{0}^{\prime }\right\rangle $ of $T_{0}^{\prime }$, 
\begin{equation*}
\rho \left\langle T_{0}^{\prime }\right\rangle =\left\langle S_{0}^{\prime
}\right\rangle .
\end{equation*}
For any such $\rho $, let $m$ be minimal such that 
\begin{equation*}
\tilde{t}^{m+1}\in \left( \rho ^{*}\mathcal{I}^{\prime }\right)
_{\left\langle T_{0}^{\prime }\right\rangle }.
\end{equation*}

We let 
\begin{equation*}
\tilde{X}/\tilde{\Delta },\ \tilde{M}/\tilde{\Delta }
\end{equation*}
be the pullbacks of $X/\Delta $ and $M/\Delta $ respectively, 
\begin{equation*}
\tilde{\Delta }_{m}=\frac{Spec\Bbb{C}\left[ \left[ \tilde{t}\right] \right] 
}{\left\{ \tilde{t}^{m+1}\right\} },
\end{equation*}
and 
\begin{equation*}
T_{m}^{\prime }=T^{\prime }\times _{\tilde{\Delta }}\tilde{\Delta }_{m}.
\end{equation*}

We let 
\begin{equation*}
\tilde{M}_{m}=\tilde{M}\times _{\tilde{\Delta }}\tilde{\Delta }_{m}.
\end{equation*}
Then, since $M_{0}$ is assumed to be the maximal abelian subvariety of$\ 
\mathcal{J}_{1}\left( X_{0}\right) $, $\tilde{M}_{m}/\tilde{\Delta }_{m}$ is
the maximal abelian subvariety of$\ \mathcal{J}_{1}\left( \tilde{X}_{m}/%
\tilde{\Delta }_{m}\right) $ for each $m$. So the Abel-Jacobi mapping 
\begin{equation*}
\tilde{\varphi}:T_{m}^{\prime }/\tilde{\Delta }_{m}\rightarrow \mathcal{J}%
_{1}\left( \tilde{X}/\tilde{\Delta }\right)
\end{equation*}
factors through the inclusion 
\begin{equation*}
\tilde{M}_{m}/\tilde{\Delta }_{m}\subseteq \mathcal{J}_{1}\left( \tilde{X}/%
\tilde{\Delta }\right) .
\end{equation*}

So we must show that, for each diagram $\left( \ref{y}\right) $ as above,
the composition 
\begin{equation}
\left( \left\{ \tilde{t}^{m+1}\right\} /\left\{ \tilde{t}^{m+2}\right\}
\right) ^{\vee }\overset{\tilde{\upsilon}}{\longrightarrow }R^{1}\tilde{p}%
_{*}N_{\left\langle T_{0}\right\rangle \backslash \left\langle T_{0}^{\prime
}\right\rangle \times X_{0}}\rightarrow \Omega _{\left\langle T_{0}^{\prime
}\right\rangle }^{1}  \label{z}
\end{equation}
is zero, that is, the map 
\begin{equation}
H^{3,0}\left( \tilde{X}/\tilde{\Delta }\right) \overset{L_{\tau }^{0,1}}{
\rightarrow }H^{2,1}\left( \tilde{X}/\tilde{\Delta }\right) \overset{\left(
restrict\right) }{\rightarrow }R^{1}\tilde{p}_{*}\left( \Omega
_{\left\langle T_{m}\right\rangle /\tilde{\Delta }_{m}}^{2}\right) =\Omega
_{\left\langle T_{m}^{\prime }\right\rangle }^{1}  \label{v}
\end{equation}
given by the Gauss-Manin connection is zero.

But in the last section we established that 
\begin{equation}
\varphi \left( \left\langle T_{0}^{\prime }\right\rangle \right) \subseteq 
\tilde{M}_{0}.  \label{claim}
\end{equation}
So by Corollary \ref{goodcor} the hypotheses of Theorem \ref{k} are
satisfied and we have established that $\left( \ref{z}\right) $ is the zero
map. Then by $\left( \ref{BBF}\right) $ we have at a general point $\left\{
Z_{0}\right\} \in S_{0}^{\prime }$ that 
\begin{equation*}
\begin{array}{r}
\dim X^{\prime }-\dim _{\left\{ Z_{0}\right\} }I^{\prime }\leq h^{1}\left(
N_{Z_{0}\backslash X_{0}}\right) -\dim S_{0}^{\prime } \\ 
=h^{0}\left( N_{Z_{0}\backslash X_{0}}\right) -\dim S_{0}^{\prime } \\ 
=\dim \left( X^{\prime }\cap \left( \left\{ 0\right\} \times J^{\prime
}\right) \right) -\dim S_{0}^{\prime }.
\end{array}
\end{equation*}
On the other hand, since codimension is lower semi-continuous, 
\begin{equation*}
\begin{array}{l}
\dim \left( X^{\prime }\cap \left( \left\{ 0\right\} \times J^{\prime
}\right) \right) -\dim S_{0}^{\prime } \\ 
=\dim \left( X^{\prime }\cap \left( \left\{ 0\right\} \times J^{\prime
}\right) \right) -\dim _{\left\{ Z_{0}\right\} }\left( I^{\prime }\cap
\left( \left\{ 0\right\} \times J^{\prime }\right) \right) \\ 
\leq \dim X^{\prime }-\dim _{\left\{ Z_{0}\right\} }I^{\prime }.
\end{array}
\end{equation*}
Thus the above inequalities are actually equalities and 
\begin{equation*}
\dim _{\left\{ Z_{0}\right\} }I^{\prime }-\dim S_{0}^{\prime }=\dim
X^{\prime }-\dim \left( X^{\prime }\cap \left( \left\{ 0\right\} \times
J^{\prime }\right) \right) =1.
\end{equation*}
So $S_{0}^{\prime }$ is a codimension-one subvariety of some component of $%
I_{red.}^{\prime }$. So, in particular, the general curve $Z_{0}$
parametrized by $S_{0}^{\prime }$ is the specialization of curves on the
nearby $X_{t}$, which completes the proof of the theorem.
\end{proof}

The weak point of Theorem \ref{g}, in addition to the fact that it is only a
variational result, is the hypothesis that the curves $Z$ of the family
parametrized by $S_{0}^{\prime }$ are unobstructed in $P$. That this
condition is necessary for the proof is shown by the following example due
to C. Voisin. Let $Z$ be the generic projection to $\Bbb{P}^{3}$ of a
canonical curve of genus 5 and let $X_{0}$ be the most general quintic in $%
\Bbb{P}^{4}$ which contains $Z$. Let $A_{0}$ be zero. In this case $%
I_{red.}^{\prime }=S_{0}^{\prime }=\Bbb{P}^{1},$ and $Z$ moves in a
basepoint-free pencil on a smooth hyperplane section of $X_{0}$. However a
constant count shows that the generic quintic threefold contains only a
finite number of curves which can specialize to curves in $I^{\prime }$. The
problem is that $Z$ lies in two components of the Hilbert scheme of curves
(of genus $5$ and degree $8$) in $\Bbb{P}^{4}$, one corresponding to full
canonical embeddings of curves of genus $5$ and the other obtained by
deforming the line bundle used to imbed $Z$ to a general line bundle of
degree $8$. (One might be tempted to reimbed $X_{0}$ ``more positively'' in
some ``bigger'' $P$ to achieve the unobstructedness of $Z$ in the new $P$,
but then, of course, one loses the vector bundle $\frak{E}$.)

Notice that the strong unobstructedness condition 
\begin{equation*}
H^{1}\left( N_{Z\backslash P}\right) =0
\end{equation*}
is always satisfied if 
\begin{equation*}
Z\subseteq P=\Bbb{P}^{n}
\end{equation*}
is a rational curve. Thus the final assertion in \ref{p}.

\section{An example}

A rather striking yet simple example of Theorem \ref{g} is the following
result on the mirror family to quintic threefolds, originally proved for the
case of lines by van Geemen \cite{AK}. Consider the pencil $X_{t}$ of
quintic threefolds given by 
\begin{equation*}
F_{t}=\left( \sum\nolimits_{j=0}^{4}x_{j}^{5}\right)
-5t\prod\nolimits_{j=0}^{4}x_{j}.
\end{equation*}
Let $\mu _{5}\leq \Bbb{C}^{*}$ denote the fifth roots-of-unity and let $G$
denote the kernel of the the group homomorphism 
\begin{eqnarray*}
\left( \mu _{5}\right) ^{\times 5} &\rightarrow &\mu _{5}. \\
\left( \xi _{0},\ldots ,\xi _{4}\right) &\longmapsto
&\prod\nolimits_{j=0}^{4}\xi _{j}
\end{eqnarray*}
Then the natural action of $G$ on each $X_{t}$ gives a family 
\begin{equation*}
Y_{t}=\frac{X_{t}}{G}
\end{equation*}
with 
\begin{eqnarray}
H^{3,0}\left( Y_{t}\right) &=&\Bbb{C\cdot }res\frac{\Omega }{F_{t}}
\label{aa} \\
H^{2,1}\left( Y_{t}\right) &=&\Bbb{C\cdot }res\frac{\left(
\prod\nolimits_{j=0}^{4}x_{j}\right) \cdot \Omega }{F_{t}^{2}}.  \notag
\end{eqnarray}
Let 
\begin{equation*}
M_{t}=\ker \left( \ \mathcal{J}_{1}\left( X_{t}\right) \rightarrow \ 
\mathcal{J}_{1}\left( Y_{t}\right) \right) .
\end{equation*}
To see that $M_{0}$ is maximal, we need only show that the $2$-dimensional
complex torus $\ \mathcal{J}_{1}\left( Y_{0}\right) $ does not contain an
elliptic curve in $H^{3,0}\left( Y_{0}\right) ^{\bot }.$ Now by $\left( \ref
{aa}\right) $ the action of the automorphism group 
\begin{equation*}
\mu _{5}=\frac{\left( \mu _{5}\right) ^{\times 5}}{G}
\end{equation*}
acts on $H^{3}\left( Y_{0}\right) $ has $1$-dimensional eigenspaces $
E_{i}=H^{4-i,i-1}\left( Y_{0}\right) $ for each of the non-trivial
characters $\xi ^{i}$ $i=1,\ldots ,4,$ of $\mu _{5}$. So the existence of
such an elliptic curve would imply that the eigenspace 
\begin{equation*}
E_{1}+E_{4},
\end{equation*}
which is defined over $\Bbb{R}$, is actually defined over $\Bbb{Q}$. This
would imply that the elliptic curve 
\begin{equation*}
\frac{H^{3,0}\left( Y_{0}\right) }{H_{3}\left( Y_{0},\Bbb{Z}\right) \cap
\left( H_{3}\left( Y_{0}\right) _{1}+H_{3}\left( Y_{0}\right) _{4}\right) }
\end{equation*}
had a non-trivial automorphism of order $5$ which is impossible.

Let $S_{0}^{\prime }$ be the Fermat quintic curve of lines on $X_{0}$ given
by 
\begin{equation*}
\alpha \left( \xi _{0},-\xi _{1},0,0,0\right) +\beta \left(
0,0,x_{2},x_{3},x_{4}\right)
\end{equation*}
such that 
\begin{equation*}
x_{2}^{5}+x_{3}^{5}+x_{4}^{5}=0
\end{equation*}
(or any image of this curve under any automorphism given by permutation of
the coordinates $x_{j}$). Thus 
\begin{equation*}
S_{0}^{\prime }\times S_{0}^{\prime }\rightarrow M_{0}\subseteq \mathcal{J}
_{1}\left( X_{0}\right) .
\end{equation*}
Thus by Theorem $\ref{g}$:

\begin{corollary}
\label{ab} Any continuous family $S_{0}^{\prime }$ of rational curves on the
Fermat quintic $X_{0}$ whose generic member is smooth lies in the
specialization at $t=0$ of a continuous family $T_{t}^{\prime }$ on $X_{t}$
for $t$ generic.
\end{corollary}

\end{document}